\documentclass[pdflatex,sn-apa]{sn-jnl}

\usepackage{graphicx}%
\usepackage{multirow}%
\usepackage{amsmath,amssymb,amsfonts}%
\usepackage{amsthm}%
\usepackage{mathrsfs}%
\usepackage[title]{appendix}%
\usepackage{xcolor}%
\usepackage{textcomp}%
\usepackage{manyfoot}%
\usepackage{booktabs}%
\usepackage{algorithm}%
\usepackage{algorithmicx}%
\usepackage{algpseudocode}%
\usepackage{listings}%

\theoremstyle{thmstyleone}%
\newtheorem{theorem}{Theorem}
\newtheorem{proposition}[theorem]{Proposition}%

\theoremstyle{thmstyletwo}%
\newtheorem{remark}{Remark}%
\newtheorem{lemma}{Lemma}

\theoremstyle{thmstylethree}%

\usepackage{float} 
\usepackage{hyperref} 
\usepackage{booktabs} 
\usepackage{booktabs} 
\usepackage{lineno} 
\usepackage{multirow} 
\usepackage{multicol} 
\usepackage{enumitem} 
\usepackage{dsfont} 
\usepackage{xcolor} 
\usepackage{appendix} 
\usepackage{pdfpages} 

\newcommand{\R}{\mathbb{R}}
\newcommand{\PP}{\mathsf{P}} 
\newcommand{\EE}{\mathsf{E}} 
\newcommand{\Bias}{\mathsf{Bias}} 
\newcommand{\Var}{\mathsf{Var}} 
\newcommand{\bb}[1]{\boldsymbol{#1}}
\newcommand{\OO}{\mathcal{O}}
\newcommand{\oo}{\mathrm{o}}
\newcommand{\rd}{\mathrm{d}}
\newcommand{\ind}{\mathds{1}}
\newcommand{\nvert}[0]{\, \vert \, }

\allowdisplaybreaks

\begin{document}

\title[Article Title]{Local linear smoothing for regression surfaces on the simplex using Dirichlet kernels}

\author*[1]{\fnm{Christian} \sur{Genest}}\email{christian.genest@mcgill.ca}

\author[1]{\fnm{Fr\'ed\'eric} \sur{Ouimet}}\email{frederic.ouimet2@mcgill.ca}

\affil*[1]{\orgdiv{Department of Mathematics and Statistics}, \orgname{McGill University}, \\
\orgaddress{\street{805, rue Sherbrooke ouest}, \city{Montr\'eal} \state{(Qu\'ebec)} \country{Canada} \postcode{H3A 0B9}}}

\abstract{This paper introduces a local linear smoother for regression surfaces on the simplex. The estimator solves a least-squares regression problem weighted by a locally adaptive Dirichlet kernel, ensuring good boundary properties. Asymptotic results for the bias, variance, mean squared error, and mean integrated squared error are derived, generalizing the univariate results of Chen [\emph{Ann. Inst. Statist. Math.}, 54(2) (2002), pp. 312--323]. A simulation study shows that the proposed local linear estimator with Dirichlet kernel outperforms its only direct competitor in the literature, the Nadaraya--Watson estimator with Dirichlet kernel due to Bouzebda, Nezzal and Elhattab [\emph{AIMS Math.}, 9(9) (2024), pp. 26195–26282].}



\maketitle

\thispagestyle{empty}


\section{Introduction}
\label{sec:1}

Although Dirichlet kernels were proposed almost $40$ years ago by \citet{doi:10.2307/2347365} for density estimation of compositional data, their application to nonparametric regression has only been considered recently. Regression in this context is valuable to predict a real response variable as a function of explanatory variables supported on the simplex (i.e., proportions that sum to 1). For instance, in Arctic lake environments, compositional data from sediment samples --- expressed as proportions of sand, silt, and clay --- vary as a function of water depth, thereby revealing underlying environmental patterns; see, e.g., \cite{Coakler/Rust:1968} and the illustration in Section~\ref{sec:6}. Another example is provided by large-scale soil composition surveys, such as the GEMAS project \citep{Reimann/etal:2014a, Reimann/etal:2014b}, which produce datasets detailing soil composition in terms of elemental concentrations; in such cases, Dirichlet kernels facilitate accurate estimation and interpretation of relationships between compositional predictors and environmental responses such as soil pH, magnetic susceptibility or annual precipitation. Similar opportunities arise in ecology, where species abundances form compositional vectors, or in health sciences, where nutrient distributions are studied. By intrinsically accommodating the simplex constraint, Dirichlet kernels ensure reliable estimation near the boundary and support a broad range of practical applications.

\citet{MR4796622} introduced, in the larger framework of conditional $U$-statistics, a version of the Nadaraya--Watson (NW) estimator weighted by a locally adaptive Dirichlet kernel, and studied some of its asymptotic properties. This development follows earlier theoretical work on Dirichlet kernel density estimators by \citet{MR4319409} and \citet{MR4544604}. \citet{MR4796622} also investigated a NW estimator weighted by a locally adaptive multinomial kernel, building on the Bernstein approach to density estimation on the simplex proposed by \citet{MR1293514} in two dimensions and expanded upon by \citet{MR4287788,doi:10.1515/stat-2022-0111} in arbitrary dimension.

As an alternative to the estimator proposed by~\citet{MR4796622}, the present paper introduces a local linear (LL) smoother for regression surfaces using Dirichlet kernels. This approach, which is introduced in Section~\ref{sec:2} along with other preliminary definitions and notations, leverages the strengths of asymmetric kernels and LL smoothing to provide an efficient estimation method on the simplex with excellent boundary properties, addressing the practical limitations of implementing the multivariate LL smoother of \citet{MR1311979} with a variable bandwidth (to our knowledge, a way to compute the variable bandwidth on the simplex is still missing).

Asymptotic results for the bias, variance, mean squared error (MSE), and mean integrated squared error (MISE) of the proposed smoother are stated in Section~\ref{sec:3} and proved in Section~\ref{sec:4}, based on technical lemmas which are relegated to the Appendix. These results, which generalize those established by \citet{MR1910175} in one dimension, are complemented in Section~\ref{sec:5} with a simulation study which demonstrates, in various settings, the superior performance of the LL estimator with Dirichlet kernel compared to its only direct competitor, the NW estimator with Dirichlet kernel recently proposed by \citet{MR4796622}. A numerical example is then presented in Section~\ref{sec:6} based on data from \cite{Coakler/Rust:1968} popularized by \cite{MR865647}.

The rest of this section provides context. First, it should be said that LL smoothers with variable bandwidth are not new; they were introduced by \citet{MR1193323} as a tool to estimate the functional relationship between a predictor and a response. Such estimators solve a locally weighted least-squares regression problem, combining the strengths of LL smoothing, as developed by \citet{MR443204} and \citet{MR556476}, with the local adaptivity resulting from the optimization of the variable bandwidth. A distinct advantage of LL estimation is that it adapts to both random and fixed designs, and exhibits no boundary effects, i.e., requires no modification near the boundary to reduce bias. As shown by \citet{MR1209561}, LL estimators also offer higher efficiency than traditional NW and Gasser--M\"uller (GM) regression estimators; their minimax properties were investigated by \citet{MR1212173}. Although LL estimation was later extended by \citet{MR1311979} to the multiple predictor setting considered herein, LL smoothers are unfortunately currently not a viable option for compositional data, due to lack of efficient methods to compute the variable bandwidth on such supports.

An alternative approach to addressing the boundary issues of traditional regression estimators was proposed by \citet{MR1742101}, who replaced the fixed kernel in GM estimators with a beta kernel. The beta kernel is locally adaptive through its parameters rather than through a variable bandwidth, i.e., the kernel parameters change with every estimation point $\bb{s}$ on the simplex. This local adaptivity allows the beta kernel to adjust its shape to fit the data within $[0,1]$, and, if appropriately transformed, within any bounded interval. \citet{MR1685301} proposed a similar approach where kernel weights are convoluted with a binomial distribution instead of partitioned, thereby offering another effective mitigation strategy for boundary bias.

The approach adopted here was inspired by the work of \citet{MR1910175}, who combined the ideas of LL smoothing and asymmetric kernels by replacing Fan's fixed kernel with variable bandwidth with a beta or gamma kernel. The former is used for design variables supported on $[0,1]$ and the latter for $[0,\infty)$. This combination ensures an intrinsic local adaptivity of the kernel and the asymptotic negligibility of boundary bias. As for Fan's LL smoother, it is suitable for both fixed and random designs.

Chen's approach to integrating locally adaptive beta/gamma kernels with LL smoothers has been influential in the literature on asymmetric kernel estimators. The results contained in the present paper extend those of \citet{MR1910175} to the multivariate case and their value is documented in the context of nonparametric regression for density estimation on the simplex.

Many other authors, including \citet{MR3256394}, \citet{MR4302589, Funke_Hirukawa_2024}, and \citet{MR4471289, MR4556820}, have expanded on Chen's ideas in different settings. When the slope in the locally weighted least-squares problem is zero, the regression estimator simplifies to a NW estimator weighted by an asymmetric kernel. This estimator has been studied by \citet{MR3494026} for gamma kernels and by \citet{MR3499725} within the broader theory of associated kernels.

\section{A local linear smoother based on Dirichlet kernels}
\label{sec:2}

Consider data $(\bb{X}_1, Y_1), \dots, (\bb{X}_n, Y_n)$, where for each $i \in \{1, \dots, n \}$, $Y_i$ is a response variable related to the random design point $\bb{X}_i$ by the functional model
\begin{equation}
\label{eq:1}
Y_i = m (\bb{X}_i) + \varepsilon_i,
\end{equation}
where $\varepsilon_i$ is a random error term and $m(\bb{x}) = \EE(Y \nvert \bb{X} = \bb{x})$ is an unknown regression function supported on the unit simplex defined, in terms of the $\ell^1$ norm $\|\cdot\|_1$ in $\R^d$, by
\[
\mathcal{S}_d = \big\{\bb{s}\in [0,1]^d : \|\bb{s}\|_1 \leq 1\big\}.
\]

Assume that the design points have a common density $f$ and that the error terms $\varepsilon_1, \dots, \varepsilon_n$ are uncorrelated with zero mean and variance $\sigma^2(\bb{X}_1), \dots, \sigma^2(\bb{X}_n)$, respectively. Further let $\sigma^2(\bb{x}) = \Var(Y \nvert \bb{X} = \bb{x})$ denote the conditional variance.

A natural local linear (LL) estimator $\hat{m}$ of $m$ can be obtained by minimizing the locally weighted loss function
\[
L(\alpha,\bb{\beta}) = \sum_{i=1}^n \{Y_i - \alpha - \bb{\beta}^{\top} (\bb{X}_i - \bb{s})\}^2 K_{\bb{s} / b + \bb{1}, (1 - \|\bb{s}\|_1) / b + 1}(\bb{X}_i),
\]
where for arbitrary $\bb{u} = (u_1, \dots, u_d) \in (0,\infty)^d$, $v \in (0,\infty)$, and $\bb{x} \in \mathcal{S}_d$,
\[
K_{\bb{u},v}(\bb{x}) = \frac{\Gamma(\|\bb{u}\|_1 + v)}{\Gamma(v) \prod_{i=1}^d \Gamma(u_i)} \, (1 - \|\bb{x}\|_1)^{v - 1} \prod_{i=1}^d x_i^{u_i - 1}
\]
is the density of the $\mathrm{Dirichlet}\hspace{0.2mm}(\bb{u},v)$ distribution.

Precisely, if $\bb{e}_1 = (1, 0, \dots, 0)$ is a $(d+1) \times 1$ vector, then one can set
\begin{equation}
\label{eq:2}
\hat{m}^{(\mathrm{LL})}_{n,b}(\bb{s}) = \hat{\alpha}_{\bb{s}} = \bb{e}_1^{\top} (\mathcal{X}_{\bb{s}}^{\top} W_{\bb{s}} \mathcal{X}_{\bb{s}})^{-1} \mathcal{X}_{\bb{s}}^{\top} W_{\bb{s}} \bb{Y},
\end{equation}
where
\[
\mathcal{X}_{\bb{s}} =
\begin{bmatrix}
1 & (\bb{X}_1 - \bb{s})^{\top} \\
\vdots & \vdots \\
1 & (\bb{X}_n - \bb{s})^{\top}
\end{bmatrix}_{n \times (d+1)}, \quad
W_{\bb{s}} = \text{diag}
\begin{bmatrix}
\kappa_{\bb{s},b}(\bb{X}_1) \\
\vdots \\
\kappa_{\bb{s},b}(\bb{X}_n)
\end{bmatrix}_{n \times 1}, \quad
\bb{Y} =
\begin{bmatrix}
Y_1 \\
\vdots \\
Y_n
\end{bmatrix}_{n \times 1},
\]
and, for simplicity, the following notation is used henceforth:
\[
\kappa_{\bb{s},b} = K_{\bb{s} / b + \bb{1}, (1 - \|\bb{s}\|_1) / b + 1}.
\]

Estimator \eqref{eq:2}, which is new, is a natural $d$-variate extension of the beta kernel estimator due to \citet{MR1910175}. As will be seen in Section~\ref{sec:3} and proved in Section~\ref{sec:4}, it enjoys excellent sampling properties. As reported in Section~\ref{sec:5}, its finite-sample behavior also makes it preferable to the alternative estimator recently proposed by \citet{MR4796622}. In addition, its practical use is illustrated in Section~\ref{sec:6} in a real-data analysis using compositional sediment samples from a Canadian Arctic lake.

In what follows, the notation $u = \OO(v)$ means that $\limsup |u / v| \leq C < \infty$ as $n\to \infty$ or $b\to 0$, depending on the context, where the positive constant $C \in (0, \infty)$ may depend on the target regression function $m$, the design density $f$, the conditional variance function $\sigma^2$, or the dimension $d$, but no other variable unless explicitly written as a subscript. Similarly, the notation $u = \oo(v)$ means that $|u / v| \to 0$ as $n \to \infty$ or $b \to 0$. Subscripts indicate on which parameters the convergence rate can depend. Except for Lemma~\ref{lem:A.b.asymptotics} in the Appendix, it is always assumed that the bandwidth parameter $b = b(n)$ depends on the sample size, $n$, in such a way that $b\to 0$ as $n\to \infty$. The shorthand $[d] = \{1,\dots,d\}$ is often used.

\section{Results}
\label{sec:3}

For any subset of indices $\mathcal{J} \subseteq [d]$, define
\begin{equation}
\label{eq:3}
\psi_{\mathcal{J}}(\bb{s}) = \left\{(4\pi)^{d - |\mathcal{J}|} (1 - \|\bb{s}\|_1) \prod_{i\in [d]\backslash\mathcal{J}} s_i\right\}^{-1/2}
\end{equation}
and
\[
h_{\mathcal{J}}(\bb{s}) =
\begin{cases}
\displaystyle \sum_{k,\ell\in [d]\backslash \mathcal{J}} \frac{1}{2} \, (s_k \ind_{\{k = \ell\}} - s_k s_{\ell}) \, \frac{\partial^2}{\partial s_k \partial s_{\ell}} \, m(\bb{s}) &\mbox{if } \mathcal{J} \subsetneq [d], \\[5mm]
\displaystyle \sum_{k,\ell\in [d]} \frac{1}{2} \{(\lambda_k + 1) \ind_{\{k = \ell\}} + 1\} \, \frac{\partial^2}{\partial s_k \partial s_{\ell}} \, m(\bb{s}) &\mbox{if } \mathcal{J} = [d], \vspace{2mm}
\end{cases}
\]
with the convention that $\prod_{\emptyset} = 1$ and $\sum_{\emptyset} = 0$. For any design point $\bb{X}$, define
\[
A_b(\bb{s}) = \EE\{\kappa_{\bb{s},b}(\bb{X})^2\}.
\]

Further assume that the following conditions hold:
\bigskip
\begin{enumerate}[label=(A\arabic*)]
\item $m\in C^2(\mathcal{S}_d)$; $f$ and $\sigma^2$ are Lipschitz continuous on $\mathcal{S}_d$. \label{ass:1} \\[-3mm]
\item $f(\bb{x}) \geq f_c > 0$ and $\sigma^2(\bb{x}) \leq \sigma_c^2 < \infty$ for all $\bb{x}\in \mathcal{S}_d$. \label{ass:2} \\[-3mm]
\item $b\to 0$ and $n b^{d/2}\to \infty$ as $n\to \infty$. \label{ass:3}
\end{enumerate}

\bigskip
Under Assumptions (A1)--(A3), Proposition~\ref{prop:1} below presents asymptotic expressions for the bias and variance of estimator \eqref{eq:2} at $\bb{s}\in \mathrm{Int}(\mathcal{S}_d)$, for the cases in which $\bb{s}$ is well inside the simplex and near the boundary. This result extends similar derivations by \citet{MR1910175} in the univariate case  $\mathcal{S}_1 \equiv [0,1]$. The index set $\mathcal{J} \subseteq [d]$ corresponds to the dimensions for which $\bb{s}$ is close the boundary of the simplex.

\bigskip
\begin{proposition}
\label{prop:1}
Suppose that Assumptions~\ref{ass:1}--\ref{ass:3} hold and let $\mathcal{J} \subseteq [d]$ and $\bb{\lambda}\in (0,\infty)^d$ be given. If a sequence $\bb{s} = \bb{s}(b)$ is selected inside $\mathrm{Int}(\mathcal{S}_d)$ such that $s_i / b \to \lambda_i$ for all $i\in \mathcal{J}$, $s_i$ is fixed for all $i\in [d]\backslash \mathcal{J}$, and $(1 - \|\bb{s}\|_1)$ is fixed, then
\[
\Bias\{\hat{m}^{(\mathrm{LL})}_{n,b}(\bb{s})\} = b^{1 + \ind_{\{\mathcal{J} = [d]\}}} h_{\mathcal{J}}(\bb{s}) + \OO\{b^{3/2 \times (1 + \ind_{\{\mathcal{J} = [d]\}})}\} + \OO\{n^{-1} A_b(\bb{s})\}
\]
and
\[
\Var\{\hat{m}^{(\mathrm{LL})}_{n,b}(\bb{s})\}
= n^{-1} b^{-(d + |\mathcal{J}|)/2} \times \left\{\frac{\psi_{\mathcal{J}}(\bb{s}) \sigma^2(\bb{s})}{f(\bb{s})} \prod_{i\in \mathcal{J}} \frac{\Gamma(2\lambda_i + 1)}{2^{2\lambda_i + 1} \Gamma^2(\lambda_i + 1)} + \oo_{\bb{\lambda},\bb{s}}(1)\right\},
\]
where the term $A_b(\bb{s})$ is bounded and estimated in Lemma~\ref{lem:A.b.asymptotics}.
\end{proposition}

\begin{remark}
In the classical setting of local polynomial smoothers with fixed kernels, it is well known that the NW estimator and the LL smoother both have $\OO(b)$ bias in the interior of the domain, but that they differ near the boundary \citep[p.~131]{MR1319818}. Specifically, in the univariate case for simplicity, the bias of the NW estimator near the boundary is $\OO(\sqrt{b})$, whereas it remains $\OO(b)$ for the LL smoother. A similar phenomenon occurs in the multivariate setting \citep{MR1311979}.

In the case of the NW estimator with Dirichlet kernel recently introduced by \citet{MR4796622}, namely
\[
\hat{m}^{(\mathrm{NW})}_{n,b}(\bb{s}) = {\sum_{i=1}^n Y_i \kappa_{\bb{s},b}(\bb{X}_i)} \Big/ {\sum_{j=1}^n \kappa_{\bb{s},b}(\bb{X}_j)},
\]
the situation differs. Because the Dirichlet kernel adapts its shape to the estimation point $\bb{s}$ and remains supported on the simplex, the bias of both the NW and LL estimators with Dirichlet kernel is at most $\OO(b)$ near the boundary.

Specifically, a detailed examination of the proof of Proposition~\ref{prop:1} given in Section~\ref{sec:4} for the NW estimator with Dirichlet kernel reveals that the main asymptotic term is driven by the ratio $E_n / A_n$ defined there. The terms $B_n$, $C_n$, $D_n$, $F_n$ also introduced in this proof do not arise in the analysis because the NW estimator is effectively a local constant estimator. One obtains
\[
\Bias\{\hat{m}^{(\mathrm{NW})}_{n,b}(\bb{s})\} = b^{1 + \ind_{\{\mathcal{J} = [d]\}}} h_{\mathcal{J}}(\bb{s}) + \OO\{n^{-1} A_b(\bb{s})\},
\]
and the additional term $\OO\{b^{3/2 \times (1 + \ind_{\{\mathcal{J} = [d]\}})}\}$ that appears in Proposition~\ref{prop:1} to account for the local slope does not arise here.

Given the superior performance of $\hat{m}^{(\mathrm{LL})}_{n,b}$ over $\hat{m}^{(\mathrm{NW})}_{n,b}$ in the simulation study in Section~\ref{sec:5}, it is conjectured that cancellations occur near the boundary between the two leading asymptotic terms of the bias in Proposition~\ref{prop:1}, i.e., in
\[
b^{1 + \ind_{\{\mathcal{J} = [d]\}}} h_{\mathcal{J}}(\bb{s}) + \OO\{b^{3/2 \times (1 + \ind_{\{\mathcal{J} = [d]\}})}\}.
\]

Tracking down an explicit closed-form expression for the higher-order asymptotics is challenging. If such cancellations do occur, they would parallel the asymptotic bias improvements observed between LL and NW estimators in the classical setting. Further investigation is required to resolve this open problem, which is not pursued further here.
\end{remark}

\bigskip
An expression for the MSE follows immediately from the above proposition. An asymptotic expression is also available for the locally optimal bandwidth under the MSE criterion. Specifically, $\smash{\mathrm{MSE}\{\hat{m}^{(\mathrm{LL})}_{n,b} (\bb{s})\}}$ is given, for any $\bb{s}\in \mathrm{Int}(\mathcal{S}_d)$, by
\begin{multline*}
n^{-1} b^{-(d + |\mathcal{J}|)/2} \, \frac{\psi_{\mathcal{J}}(\bb{s}) \sigma^2(\bb{s})}{f(\bb{s})} \prod_{i\in \mathcal{J}} \frac{\Gamma(2\lambda_i + 1)}{2^{2\lambda_i + 1} \Gamma^2(\lambda_i + 1)} + b^{2 + 2\ind_{\{\mathcal{J} = [d]\}}} \, \{h_{\mathcal{J}}(\bb{s})\}^2 \\[1mm]
+ \OO_{\bb{\lambda},\bb{s}}(n^{-1} b^{-(d + |\mathcal{J}| - 1)/2}) + \OO\{b^{5/2 (1 + \ind_{\{\mathcal{J} = [d]\}})}\} \\[2mm]
+ \OO\{b^{1 + \ind_{\{\mathcal{J} = [d]\}}} n^{-1} A_b(\bb{s})\} + \OO\{n^{-2} A_b(\bb{s})^2\}.
\end{multline*}

In particular, if $h_{\mathcal{J}}(\bb{s}) \neq 0$ and $R_{\mathcal{J}} \equiv d + |\mathcal{J}| + 4 + 4 \ind_{\{\mathcal{J} = [d]\}}$, the asymptotically optimal choice of $b$ with respect to $\mathrm{MSE}$ is
\[
b_{\mathrm{opt}}(\bb{s}) = \Bigg[n^{-1} \times \frac{d + |\mathcal{J}|}{4 + 4 \ind_{\{\mathcal{J} = [d]\}}} \times \frac{\psi_{\mathcal{J}}(\bb{s}) \sigma^2(\bb{s})}{f(\bb{s}) \{h_{\mathcal{J}}(\bb{s})\}^2} \prod_{i\in \mathcal{J}} \frac{\Gamma(2\lambda_i + 1)}{2^{2\lambda_i + 1} \Gamma^2(\lambda_i + 1)}\Bigg]^{2/R_{\mathcal{J}}},
\]
and the corresponding value of the MSE is
\begin{multline*}
\mathrm{MSE}\{\hat{m}^{(\mathrm{LL})}_{n,b_{\mathrm{opt}}(\bb{s})}(\bb{s})\} \\
= n^{-(4 + 4 \ind_{\{\mathcal{J} = [d]\}})/R_{\mathcal{J}}} \left(1 + \frac{d + |\mathcal{J}|}{4 + 4 \ind_{\{\mathcal{J} = [d]\}}}\right) \left(\frac{d + |\mathcal{J}|}{4 + 4 \ind_{\{\mathcal{J} = [d]\}}}\right)^{-(d + |\mathcal{J}|)/R_{\mathcal{J}}} \\
\qquad \qquad \times \left\{\frac{\psi_{\mathcal{J}}(\bb{s}) \sigma^2(\bb{s})}{f(\bb{s})} \prod_{i\in \mathcal{J}} \frac{\Gamma(2\lambda_i + 1)}{2^{2\lambda_i + 1} \Gamma^2(\lambda_i + 1)}\right\}^{(4 + 4 \ind_{\{\mathcal{J} = [d]\}})/R_{\mathcal{J}}} \\
\times \big[\{h_{\mathcal{J}}(\bb{s})\}^2\big]^{(d + |\mathcal{J}|)/R_{\mathcal{J}}} + \oo_{\bb{s}}\Big(n^{-(4 + 4 \ind_{\{\mathcal{J} = [d]\}})/R_{\mathcal{J}}}\Big).
\end{multline*}

By integrating the MSE and showing that the contribution from points near the boundary is negligible, one obtains an expression for the MISE, viz.
\[
\mathrm{MISE}\{\hat{m}^{(\mathrm{LL})}_{n,b}\} = \int_{\mathcal{S}_d} \EE\big\{|\hat{m}^{(\mathrm{LL})}_{n,b}(\bb{s}) - m(\bb{s})|^2\big\} \, \rd \bb{s}.
\]
This expression, detailed below, could be used to implement a plug-in method; see, e.g., Section~3.6 of~\citet{MR3822372}.

\bigskip
\begin{proposition}
\label{prop:2}
Under Assumptions~\ref{ass:1}--\ref{ass:3}, one has
\[
\mathrm{MISE}\{\hat{m}^{(\mathrm{LL})}_{n,b}\}
= \frac{1}{n b^{d/2}} \int_{\mathcal{S}_d} \frac{\psi_{\emptyset}(\bb{s}) \sigma^2(\bb{s})}{f(\bb{s})} \, \rd \bb{s} + b^2 \int_{\mathcal{S}_d} \{h_{\emptyset}(\bb{s})\}^2 \rd \bb{s} + \oo\left(\frac{1}{n b^{d/2}}\right) + \oo(b^2).
\]
In particular, if $\int_{\mathcal{S}_d} \{h_{\emptyset}(\bb{s})\}^2 \rd \bb{s} > 0$, the asymptotically optimal choice of $b$ with respect to $\mathrm{MISE}$, denoted $b_{\mathrm{opt}}(\bb{s})$, is given by
\begin{multline*}
 n^{-2/(d+4)} (d/4)^{2/(d+4)} \left\{\int_{\mathcal{S}_d} \frac{\psi_{\emptyset}(\bb{s}) \sigma^2(\bb{s})}{f(\bb{s})} \, \rd \bb{s}\right\}^{2/(d+4)} \left[\int_{\mathcal{S}_d} \{h_{\emptyset}(\bb{s})\}^2 \rd \bb{s}\right]^{-2/(d+4)},
\end{multline*}
and the corresponding value of the MISE is
\begin{multline*}
\mathrm{MISE}\{\hat{m}^{(\mathrm{LL})}_{n,b_{\mathrm{opt}}}\}
= n^{-4 / (d+4)} \left\{\frac{1 + d/4}{(d/4)^{d/(d+4)}}\right\} \left\{\int_{\mathcal{S}_d} \frac{\psi_{\emptyset}(\bb{s}) \sigma^2(\bb{s})}{f(\bb{s})} \, \rd \bb{s}\right\}^{4/(d+4)} \\ \times \left[\int_{\mathcal{S}_d} \{h_{\emptyset}(\bb{s})\}^2 \rd \bb{s}\right]^{d/(d+4)} + \oo\big(n^{-4/(d+4)}\big).
\end{multline*}
\end{proposition}

\section{Mathematical arguments}
\label{sec:4}

This section contains detailed proofs of Propositions~\ref{prop:1} and \ref{prop:2}.

\subsection{Proof of Proposition~\ref{prop:1}}
\label{sec:4.1}

To prove Proposition~\ref{prop:1}, first note that Assumption~\ref{ass:1} and a stochastic Taylor expansion (see Theorem~18.18 of \citet{MR2378491}) together imply that
\[
\EE \big( \bb{Y} \nvert \bb{X}_1, \dots, \bb{X}_n \big)
=
\begin{bmatrix}
m(\bb{X}_1) \\
\vdots \\
m(\bb{X}_n)
\end{bmatrix}_{n \times 1}
= \,
\mathcal{X}_{\bb{s}}
\begin{bmatrix}
m(\bb{s}) \\[1mm]
\frac{\partial}{\partial s_1} \, m(\bb{s}) \\
\vdots \\
\frac{\partial}{\partial s_d} \, m(\bb{s})
\end{bmatrix}_{(d+1) \times 1} + \, \tfrac{1}{2} \, \Lambda_m(\bb{s}),
\]
where
\[
 \Lambda_m(\bb{s}) =
\left[\hspace{-1mm}
\begin{array}{c}
\Lambda_{m,1}(\bb{s}) \\
\vdots \\
\Lambda_{m,n}(\bb{s})
\end{array}
\hspace{-1mm}\right]_{n\times 1}
\]
and
\[
\Lambda_{m, i}(\bb{s}) = \sum_{k,\ell=1}^d (X_{i,k} - s_k) (X_{i,\ell} - s_{\ell}) \, \frac{\partial^2}{\partial s_k \partial s_{\ell}} \, m(\bb{s}) \{1 + \oo_{\PP}(1)\}.
\]

Using the definition of the estimator in \eqref{eq:2}, one can write
\begin{equation}
\label{eq:4}
\EE\big\{\hat{m}^{(\mathrm{LL})}_{n,b}(\bb{s}) - m(\bb{s}) \nvert \bb{X}_1,\dots, \bb{X}_n\big\} = \frac{1}{2} \, \bb{e}_1^{\top} (n^{-1} \mathcal{X}_{\bb{s}}^{\top} W_{\bb{s}} \mathcal{X}_{\bb{s}})^{-1} n^{-1} \mathcal{X}_{\bb{s}}^{\top} W_{\bb{s}} \Lambda_m(\bb{s}).
\end{equation}

If $\bb{\xi}_{\bb{s}}\sim \mathrm{Dirichlet}\hspace{0.2mm}[\bb{s} / b + \bb{1}, (1 - \|\bb{s}\|_1) / b + 1]$, then Eq.~(4) in~\citet{MR4319409} shows that
\begin{equation}
\label{eq:5}
\EE(\bb{\xi}_{\bb{s}}) = \bb{s} + \OO(b), \quad \EE\{(\bb{\xi}_{\bb{s}} - \bb{s}) (\bb{\xi}_{\bb{s}} - \bb{s})^{\top}\} = b \{\mathrm{diag}(\bb{s}) - \bb{s} \bb{s}^{\top}\} + \OO(b^2).
\end{equation}

By Assumption~\ref{ass:1}, the weak law of large numbers and Chebyshev's inequality to control the probabilistic errors, it follows that
\[
(n^{-1} \mathcal{X}_{\bb{s}}^{\top} W_{\bb{s}} \mathcal{X}_{\bb{s}})^{-1} = \left[ \begin{array}{cc} A_n & B_n \\[0mm] C_n & D_n \end{array} \right]^{-1},
\]
where
\[
A_n = \frac{1}{n} \sum_{i=1}^n \kappa_{\bb{s},b}(\bb{X}_i) = f(\bb{s}) + \OO_{\PP}\{n^{-1} A_b(\bb{s})\} + \sum_{k=1}^d \OO\{\EE(|\xi_{\bb{s},k} - s_k|)\},
\]
\begin{multline*}
B_n = \frac{1}{n} \sum_{i=1}^n (\bb{X}_i - \bb{s})^{\top} \kappa_{\bb{s},b}(\bb{X}_i) \\
= \EE(\bb{\xi}_{\bb{s}} - \bb{s})^{\top} f(\bb{s}) + \oo_{\PP}\{n^{-1} A_b(\bb{s})\}
+ \sum_{k=1}^d \OO\{\EE(|\bb{\xi}_{\bb{s}} - \bb{s}|^{\top} |\xi_{\bb{s},k} - s_k|)\},
\end{multline*}
\begin{multline*}
C_n = \frac{1}{n}\sum_{i=1}^n (\bb{X}_i - \bb{s}) \kappa_{\bb{s},b}(\bb{X}_i) \\
= \EE(\bb{\xi}_{\bb{s}} - \bb{s}) \, f(\bb{s}) + \oo_{\PP}\{n^{-1} A_b(\bb{s})\}
+ \sum_{k=1}^d \OO\{\EE(|\bb{\xi}_{\bb{s}} - \bb{s}| |\xi_{\bb{s},k} - s_k|\},
\end{multline*}
and
\begin{multline*}
D_n = \frac{1}{n} \sum_{i=1}^n (\bb{X}_i - \bb{s}) (\bb{X}_i - \bb{s})^{\top} \kappa_{\bb{s},b}(\bb{X}_i) \\
= \EE\{(\bb{\xi}_{\bb{s}} - \bb{s}) (\bb{\xi}_{\bb{s}} - \bb{s})^{\top}\} \, f(\bb{s}) + \oo_{\PP}\{n^{-1} A_b(\bb{s})\} \\
+ \sum_{k=1}^d \OO\{\EE(|\bb{\xi}_{\bb{s}} - \bb{s}| |\bb{\xi}_{\bb{s}} - \bb{s}|^{\top} |\xi_{\bb{s},k} - s_k|)\}.
\end{multline*}

Therefore, using \eqref{eq:5}, and then Lemma~\ref{lem:inverse.block.matrix} to invert the block matrix, one obtains
\begin{multline}
\label{eq:6}
(n^{-1} \mathcal{X}_{\bb{s}}^{\top} W_{\bb{s}} \mathcal{X}_{\bb{s}})^{-1} \\
= \left[ \begin{array}{c|c}
f(\bb{s})^{-1} + \OO_{\PP}\{n^{-1} A_b(\bb{s})\} & \OO(b \bb{1}^{\top}) + \oo_{\PP}\{n^{-1} A_b(\bb{s}) \, \bb{1}^{\top}\} \\[1.5mm]
\midrule
\OO(b \bb{1}) + \oo_{\PP}\{n^{-1} A_b(\bb{s}) \, \bb{1}\} & \star
\end{array}
\right]_{(d+1)\times (d+1)},
\end{multline}
where the expression replacing the symbol $\star$ is irrelevant. It should be noted in these calculations that the inverse covariance matrix, $\{\mathrm{diag}(\bb{s}) - \bb{s} \bb{s}^{\top}\}^{-1}$, is equal to $\mathrm{diag}(s_1^{-1},\ldots,s_d^{-1}) + s_{d+1}^{-1} \bb{1} \bb{1}^{\top}$ by Eq.~(21) of \citet{MR1157720}, and that it is bounded away from zero in each component regardless of the choice of sequence $\bb{s} = \bb{s}(b)$ in the statement of the proposition.

Moreover, one has
\begin{equation}
\label{eq:7}
n^{-1} \mathcal{X}_{\bb{s}}^{\top} W_{\bb{s}} \Lambda_m(\bb{s}) = \left[ \begin{array}{cc} E_n \\ F_n \end{array} \hspace{-1mm}\right]_{(d+1)\times 1},
\end{equation}
where
\[
E_n = \frac{1}{n} \sum_{i=1}^n \sum_{k,\ell=1}^d (X_{i,k} - s_k) (X_{i,\ell} - s_{\ell}) \kappa_{\bb{s},b}(\bb{X}_i) \frac{\partial^2}{\partial s_k \partial s_{\ell}} \, m(\bb{s}) \{1 + \oo_{\PP}(1)\},
\]
and
\begin{multline*}
F_n = \frac{1}{n} \sum_{i=1}^n \sum_{k,\ell=1}^d (X_{i,k} - s_k) (X_{i,\ell} - s_{\ell}) (\bb{X}_i - \bb{s}) \kappa_{\bb{s},b}(\bb{X}_i) \frac{\partial^2}{\partial s_k \partial s_{\ell}} \, m(\bb{s}) \{1 + \oo_{\PP}(1)\}.
\end{multline*}

In view of \eqref{eq:6} and \eqref{eq:7}, one can then re-express \eqref{eq:4} as
\begin{multline*}
 \frac{f(\bb{s})^{-1}}{2n} \sum_{k,\ell=1}^d \sum_{i=1}^n (X_{i,k} - s_k) (X_{i,\ell} - s_{\ell}) \kappa_{\bb{s},b}(\bb{X}_i) \frac{\partial^2}{\partial s_k \partial s_{\ell}} \, m(\bb{s}) \{1 + \oo_{\PP}(1)\} \\
\quad+ b \sum_{j,k,\ell=1}^d \OO\left\{\frac{1}{n} \sum_{i=1}^n |X_{i,j} - s_j| |X_{i,k} - s_k| |X_{i,\ell} - s_{\ell}| \kappa_{\bb{s},b}(\bb{X}_i)\right\} + \OO_{\PP}\{n^{-1} A_b(\bb{s})\}.
\end{multline*}

Now, H\"older's inequality implies that, for all $j, k, \ell \in [d]$,
\begin{multline*}
\EE(|\xi_{\bb{s},j} - s_j| |\xi_{\bb{s},k} - s_k| |\xi_{\bb{s},\ell} - s_{\ell}|) \leq \max_{j \in [d]} \{\EE(|\xi_{\bb{s},j} - s_j|^6)\}^{3/6} = \OO\{(b^3)^{3/6}\} = \OO(b^{3/2}).
\end{multline*}

Therefore, by Assumption~\ref{ass:1}, the weak law of large numbers and Chebyshev's inequality to control the probabilistic errors, taking the expectation in \eqref{eq:4} yields
\begin{multline*}
\EE\big\{\hat{m}^{(\mathrm{LL})}_{n,b}(\bb{s}) - m(\bb{s})\big\}
= \EE\big[\EE\big\{\hat{m}^{(\mathrm{LL})}_{n,b}(\bb{s}) - m(\bb{s}) \nvert \bb{X}_1, \dots, \bb{X}_n\big\}\big] \\
= \frac{1}{2} \sum_{k,\ell=1}^d \EE\{(\xi_{\bb{s},k} - s_k) (\xi_{\bb{s},\ell} - s_{\ell})\} \frac{\partial^2}{\partial s_k \partial s_{\ell}} \, m(\bb{s}) \\
+ \sum_{j,k,\ell=1}^d \OO\{\EE(|\xi_{\bb{s},j} - s_j| |\xi_{\bb{s},k} - s_k| |\xi_{\bb{s},\ell} - s_{\ell}|)\} + \OO\{n^{-1} A_b(\bb{s})\},
\end{multline*}
which reduces to
\begin{align}
\label{eq:8}
\frac{b}{2} \sum_{k,\ell=1}^d (s_k \ind_{\{k = \ell\}} - s_k s_{\ell}) \frac{\partial^2}{\partial s_k \partial s_{\ell}} \, m(\bb{s}) + \OO(b^{3/2}) + \OO\{n^{-1} A_b(\bb{s})\}.
\end{align}

This establishes the expression for the bias when $\mathcal{J} = \emptyset$. More generally, note that
\[
\EE\{(\xi_{\bb{s},k} - s_k) (\xi_{\bb{s},\ell} - s_{\ell})\} = \EE(\xi_{\bb{s},k} \xi_{\bb{s},\ell}) - s_k \EE(\xi_{\bb{s},\ell}) - s_{\ell} \EE(\xi_{\bb{s},k}) + s_k s_{\ell}
\]
can be expressed as
\begin{multline*}
\frac{(s_k / b + 1) \ind_{\{k = \ell\}} - \frac{(s_k / b + 1) (s_{\ell} / b + 1)}{(1 / b + d + 1)}}{(1 / b + d + 1) (1 / b + d + 2)} + \frac{s_k / b + 1}{1 / b + d + 1} \times \frac{s_{\ell} / b + 1}{1 / b + d + 1} \\
- s_k \frac{s_{\ell} / b + 1}{1 / b + d + 1} - s_{\ell} \frac{s_k / b + 1}{1 / b + d + 1} + s_k s_{\ell}.
\end{multline*}

The latter expression, say $C_b$, is such that
\[
\begin{aligned}
C_b = \begin{cases}
b \, (s_k \ind_{\{k = \ell\}} - s_ i s_{\ell}) + \OO(b^2) &\mbox{if } s_k ~\text{is fixed and } s_{\ell} ~\text{is fixed}, \\
b^2 \{1 - s_{\ell} (\lambda_k + d + 2)\} + \OO(b^3) &\mbox{if } s_k / b \to \lambda_k ~\text{and}~ s_{\ell} ~\text{is fixed}, \\
b^2 \{1 - s_k (\lambda_{\ell} + d + 2)\} + \OO(b^3) &\mbox{if } s_k ~\text{is fixed and } s_{\ell} / b \to \lambda_{\ell}, \\
b^2 \{(\lambda_k + 1) \ind_{\{k = \ell\}} + 1\} + \OO(b^3) &\mbox{if } s_k / b \to \lambda_k ~\text{and}~ s_{\ell} / b \to \lambda_{\ell}.
\end{cases}
\end{aligned}
\]

Also, if $\mathcal{J} = [d]$, then H\"older's inequality implies that, for all $j, k, \ell \in [d]$,
\[
\EE(|\xi_{\bb{s},j} - s_j| |\xi_{\bb{s},k} - s_k| |\xi_{\bb{s},\ell} - s_{\ell}|) \leq \max_{j \in [d]} \{\EE(|\xi_{\bb{s},j} - s_j|^6)\}^{3/6} = \OO\{(b^6)^{3/6}\} = \OO(b^3).
\]
Hence, for any $\mathcal{J}\subseteq [d]$, one can conclude from \eqref{eq:8} that
\[
\EE \big\{\hat{m}^{(\mathrm{LL})}_{n,b}(\bb{s}) - m(\bb{s}) \big\} = b^{1 + \ind_{\{\mathcal{J} = [d]\}}} h_{\mathcal{J}}(\bb{s}) + \OO\{b^{3/2 \times (1 + \ind_{\{\mathcal{J} = [d]\}})}\} + \OO\{n^{-1} A_b(\bb{s})\}.
\]

To estimate the conditional variance of $\hat{m}^{(\mathrm{LL})}_{n,b}(\bb{s})$, let
\[
V = \text{diag}\{\sigma^2(\bb{X}_1),\dots,\sigma^2(\bb{X}_n)\}.
\]
By \eqref{eq:2}, one has
\begin{multline}
\label{eq:9}
\Var \big\{\hat{m}^{(\mathrm{LL})}_{n,b}(\bb{s}) \nvert \bb{X}_1, \dots, \bb{X}_n \big\} \\
= n^{-1} \bb{e}_1^{\top} (n^{-1} \mathcal{X}_{\bb{s}}^{\top} W_{\bb{s}} \mathcal{X}_{\bb{s}})^{-1} n^{-1} \mathcal{X}_{\bb{s}}^{\top} W_{\bb{s}} V W_{\bb{s}} \mathcal{X}_{\bb{s}} (n^{-1} \mathcal{X}_{\bb{s}}^{\top} W_{\bb{s}} \mathcal{X}_{\bb{s}})^{-1} \bb{e}_1.
\end{multline}

Furthermore, one can write
\[
n^{-1} \mathcal{X}_{\bb{s}}^{\top} W_{\bb{s}} V W_{\bb{s}} \mathcal{X}_{\bb{s}} = \left[ \begin{array}{cc} G_n & H_n \\ I_n & J_n \end{array} \right], \vspace{-2mm}
\]
where
\[
\begin{aligned}
G_n = \frac{1}{n} \sum_{i=1}^n \kappa_{\bb{s},b}(\bb{X}_i)^2 \, \sigma^2(\bb{X}_i),
&~~H_n = \frac{1}{n} \sum_{i=1}^n (\bb{X}_i - \bb{s})^{\top} \kappa_{\bb{s},b}(\bb{X}_i)^2 \, \sigma^2(\bb{X}_i), \\
I_n = \frac{1}{n} \sum_{i=1}^n (\bb{X}_i - \bb{s}) \kappa_{\bb{s},b}(\bb{X}_i)^2 \, \sigma^2(\bb{X}_i),
&~~J_n = \frac{1}{n} \sum_{i=1}^n (\bb{X}_i - \bb{s}) (\bb{X}_i - \bb{s})^{\top} \kappa_{\bb{s},b}(\bb{X}_i)^2 \, \sigma^2(\bb{X}_i),
\end{aligned}
\]
and so
\[
n^{-1} \mathcal{X}_{\bb{s}}^{\top} W_{\bb{s}} V W_{\bb{s}} \mathcal{X}_{\bb{s}} =
\left[\hspace{-1mm}
\begin{array}{c|c}
\begin{array}{l}
n^{-1} A_b(\bb{s}) \sigma^2(\bb{s}) f(\bb{s}) + \oo_{\PP}\{n^{-1} A_b(\bb{s})\}
\end{array}
&
\begin{array}{l}
\oo_{\PP}\{n^{-1} A_b(\bb{s})\}
\end{array} \\[1.5mm]
\midrule
\begin{array}{l}
\oo_{\PP}\{n^{-1} A_b(\bb{s})\}
\end{array}
&
\begin{array}{l}
\star
\end{array}
\end{array}
\hspace{-1mm}\right],
\]
where the expression replacing the symbol $\star$ is again irrelevant. Under Assumption~\ref{ass:3}, note that, for any given $\bb{s}\in \mathrm{Int}(\mathcal{S}_d)$, one has $A_b(\bb{s})/n \to 0$ as $n \to \infty$. Together with \eqref{eq:6} and \eqref{eq:9}, one gets
\[
\Var\big\{\hat{m}^{(\mathrm{LL})}_{n,b}(\bb{s}) \nvert \bb{X}_1, \dots, \bb{X}_n\big\} = n^{-1} {A_b(\bb{s}) \sigma^2(\bb{s})} / {f(\bb{s})} + \oo_{\PP}\{n^{-1} A_b(\bb{s})\}.
\]

It is also possible to show that
\[
\Var\big[\EE\big\{\hat{m}^{(\mathrm{LL})}_{n,b}(\bb{s}) \nvert \bb{X}_1, \dots, \bb{X}_n\big\}\big] = \oo\{n^{-1} A_b(\bb{s}) \}.
\]

Given that
\[
\Var\big\{\hat{m}^{(\mathrm{LL})}_{n,b}(\bb{s})\big\} = \Var\big[\EE\big\{\hat{m}^{(\mathrm{LL})}_{n,b}(\bb{s}) \nvert \bb{X}_1, \dots, \bb{X}_n\big\}\big] + \EE\big[\Var\big\{\hat{m}^{(\mathrm{LL})}_{n,b}(\bb{s}) \nvert \bb{X}_1, \dots, \bb{X}_n\big\}\big],
\]
it follows that
\begin{equation}
\label{eq:10}
\Var\big\{\hat{m}^{(\mathrm{LL})}_{n,b}(\bb{s})\big\} = n^{-1} {A_b(\bb{s}) \sigma^2(\bb{s})} / {f(\bb{s})} + \oo\{n^{-1} A_b(\bb{s})\}.
\end{equation}

The claimed expression for the variance follows from Lemma~\ref{lem:A.b.asymptotics} in the Appendix. This concludes the proof of Proposition~\ref{prop:1}.

\subsection{Proof of Proposition~\ref{prop:2}}
\label{sec:4.2}

Turning to the proof of Proposition~\ref{prop:2}, one can apply the bound and convergence of $A_b$ in Lemma~\ref{lem:A.b.asymptotics}, invoke the boundedness assumptions on $f$ and $\sigma^2$ from Assumption~\ref{ass:2}, and use the dominated convergence theorem to deduce from \eqref{eq:10} that
\begin{align*}
n b^{\hspace{0.2mm}d/2} \int_{\mathcal{S}_d} \Var\big\{\hat{m}^{(\mathrm{LL})}_{n,b}(\bb{s})\big\} \rd \bb{s} & = b^{\hspace{0.2mm}d/2} \int_{\mathcal{S}_d} \frac{A_b(\bb{s}) \sigma^2(\bb{s})}{f(\bb{s})} \, \{1 + \oo(1)\} \, \rd \bb{s} \\
& = \int_{\mathcal{S}_d} \frac{\psi_{\emptyset}(\bb{s}) \sigma^2(\bb{s})}{f(\bb{s})} \,\rd \bb{s} + \oo(1).
\end{align*}

Similarly, noting that $m\in C^2(\mathcal{S}_d)$ by Assumption~\ref{ass:1}, and $\mathcal{S}_d$ is compact, one can use the boundedness of the second order partial derivatives of $m$, together with the dominated convergence theorem, to deduce from Proposition~\ref{prop:1} that
\[
b^{-2} \int_{\mathcal{S}_d} \Bias\{\hat{m}_{n,b}(\bb{s})\}^2 \, \rd \bb{s} = \int_{\mathcal{S}_d} h_{\emptyset}^2(\bb{s}) \, \rd \bb{s} + \oo(1) + \OO\{b^{-2} (n^{-1} b^{-d/2})^2\}.
\]

By combining the last two equations, one finds that
\begin{align*}
\mathrm{MISE}\{\hat{m}^{(\mathrm{LL})}_{n,b}\}
& = \int_{\mathcal{S}_d} \Var\{\hat{m}_{n,b}(\bb{s})\} \rd \bb{s} + \int_{\mathcal{S}_d} \Bias\{\hat{m}_{n,b}(\bb{s})\}^2 \rd \bb{s} \\
& = n^{-1} b^{-d/2} \int_{\mathcal{S}_d} \frac{\psi_{\emptyset}(\bb{s}) \sigma^2(\bb{s})}{f(\bb{s})} \, \rd \bb{s} + b^2 \int_{\mathcal{S}_d} h_{\emptyset}^2(\bb{s}) \rd \bb{s} + \oo(b^2) + \oo(n^{-1} b^{-d/2}),
\end{align*}
where the last error term results from the fact that $(n^{-1} b^{-d/2})^2 = \oo(n^{-1} b^{-d/2})$ by Assumption~\ref{ass:3}. This concludes the proof of Proposition~\ref{prop:2}.

\section{Simulation study}
\label{sec:5}

This section reports the results of a simulation study comparing the performance of the LL estimator with Dirichlet kernel, $\smash{\hat{m}^{(\mathrm{LL})}_{n,b}}$, against its only direct competitor in the literature, the Nadaraya--Watson (NW) estimator with Dirichlet kernel introduced recently by \citet{MR4796622}, viz.
\[
\hat{m}^{(\mathrm{NW})}_{n,b}(\bb{s}) = {\sum_{i=1}^n Y_i \kappa_{\bb{s},b}(\bb{X}_i)} {\Big /} {\sum_{j=1}^n \kappa_{\bb{s},b}(\bb{X}_j)}.
\]

The bandwidth selection procedure considered is least-squares cross validation (LSCV). For a given method in the set $\{\mathrm{LL},\mathrm{NW}\}$, and a given target regression function $m$, the bandwidth $\hat{b}_n\in (0,\infty)$ is chosen to minimize the criterion
\[
\mathrm{LSCV}(b) = \frac{1}{1000} \sum_{i=1}^{1000} |\hat{m}_{n,b}^{(\mathrm{method})}(\bb{U}_i) - m(\bb{U}_i)|^2 / d!,
\]
where $\bb{U}_1,\ldots,\bb{U}_{1000}$ form a random sample from the uniform distribution $\mathcal{U}(\mathcal{S}_d)$ on the unit simplex, $\mathcal{S}_d$. The factor $1/d!$ in the numerator is the normalization constant for this distribution. For each $b \in (0,\infty)$, $\mathrm{LSCV}(b)$ is an unbiased estimator of
\[
\mathrm{MISE}\{\hat{m}_{n,b}^{(\mathrm{method})}\} \equiv \smash{\int_{\mathcal{S}_d} \EE\big\{|\hat{m}_{n,b}^{(\mathrm{method})}(\bb{s}) - m(\bb{s})|^2\big\} \rd \bb{s}}.
\]

Six target regression functions were tested, respectively defined, for all $\bb{s} = (s_1,s_2)\in \mathcal{S}_2$, by
\begin{multicols}{2}
\begin{itemize}\setlength\itemsep{0.3em}
    \item [(i)] $m_1(\bb{s}) = s_1 s_2$;
    \item [(ii)] $m_2(\bb{s}) = \ln (1 + s_1 + s_2)$;
    \item [(iii)] $m_3(\bb{s}) = \sin(s_1) + \cos(s_2)$;
    \item [(iv)] $m_4(\bb{s}) = \sqrt{s_1} + \sqrt{s_2}$;
    \item [(v)] $m_5(\bb{s}) = (s_1 + 1/4)^2 + (s_2 + 3/4)^2$;
    \item [(vi)] $m_6(\bb{s}) = (1 + s_1) e^{s_2}$.
\end{itemize}
\end{multicols}
The target regression functions were selected to represent a range of behaviors both inside the simplex and near its boundary, while conforming to the smoothness assumptions made in Section~\ref{sec:3}. They include variations in functional form (polynomial, logarithmic, trigonometric, square-root, and exponential), curvature, growth rates, boundary slopes, symmetry, and overall complexity.

Let $R = 100$ be the number of Monte Carlo replications, and suppose that the sequence of design points $\bb{x}_1,\ldots,\bb{x}_n$ consists of the fixed mesh
\[
\mathcal{M}_k = \{(w_k (i - 1) + 1/2, w_k (k - j) + 1/2)/(k+1) : 1\leq i \leq j \leq k \}
\]
inside the simplex, for $w_k = (k - 1 / \sqrt{2}) / (k - 1)$ and some positive integer $k$. For each method in $\{\mathrm{LL},\mathrm{NW}\}$, each target regression function $m_j$,~$j\in \{1,\ldots,6\}$, and each sample size $n = k (k+1)/2$ with $k\in \{7, 10, 14, 20\}$, the mean, median, standard deviation, and interquartile range, are reported in Table~\ref{tab:1} for the sequence of Monte Carlo estimates of the integrated squared errors, viz.
\[
\widetilde{\mathrm{ISE}}_{\mathrm{method}}^{j,r} = \frac{1}{1000} \sum_{i=1}^{1000} |\hat{m}_{n,\hat{b}_{n,r}}^{(\mathrm{method})}(\bb{U}_i) - m_j(\bb{U}_i)|^2 / 2,
\]
where $r \in \{1,\ldots,R\}$. In each case $(m_j,n$), Table~\ref{tab:1} shows that the new LL estimator with Dirichlet kernel has the smallest mean, median, standard deviation, and interquartile range, clearly demonstrating its dominance over the NW estimator with Dirichlet kernel of \citet{MR4796622}.

\begin{table}[ht]
\caption{Comparison of the LL and NW methods, using the fixed design points $\bb{x}_1,\ldots,\bb{x}_n$, in terms of the mean, median, standard deviation (SD), and interquartile range (IQR) of $100$ $\smash{\widetilde{\mathrm{ISE}}}$ values for regression functions $m_1, \ldots, m_6$ and sample size $n \in \{28, 55, 105, 210\}$.}
\label{tab:1}

\bigskip
\renewcommand{\arraystretch}{1.05} 
\setlength{\tabcolsep}{3.3pt} 
\centering
\begin{tabular}{cc|cccccccc}
\toprule
\multirow{2}{*}{Function} & \multirow{2}{*}{$n$} & \multicolumn{2}{c}{Mean} & \multicolumn{2}{c}{Median} & \multicolumn{2}{c}{SD} & \multicolumn{2}{c}{IQR} \\
\cmidrule(lr){3-4} \cmidrule(lr){5-6} \cmidrule(lr){7-8} \cmidrule(lr){9-10}
 & & LL & NW & LL & NW & LL & NW & LL & NW \\
\midrule
$m_1$ & \phantom{1}28 & 0.000030 & 0.000106 & 0.000030 & 0.000106 & 0.000002 & 0.000006 & 0.000002 & 0.000008 \\
$m_1$ & \phantom{1}55 & 0.000020 & 0.000061 & 0.000020 & 0.000061 & 0.000001 & 0.000003 & 0.000002 & 0.000005 \\
$m_1$ & 105 & 0.000014 & 0.000031 & 0.000014 & 0.000031 & 0.000001 & 0.000002 & 0.000001 & 0.000003 \\
$m_1$ & 210 & 0.000011 & 0.000018 & 0.000011 & 0.000018 & 0.000001 & 0.000001 & 0.000001 & 0.000002 \\
\midrule
$m_2$ & \phantom{1}28 & 0.000040 & 0.000231 & 0.000040 & 0.000229 & 0.000002 & 0.000015 & 0.000002 & 0.000024 \\
$m_2$ & \phantom{1}55 & 0.000035 & 0.000183 & 0.000035 & 0.000183 & 0.000001 & 0.000010 & 0.000001 & 0.000015 \\
$m_2$ & 105 & 0.000024 & 0.000124 & 0.000024 & 0.000124 & 0.000001 & 0.000006 & 0.000001 & 0.000009 \\
$m_2$ & 210 & 0.000014 & 0.000065 & 0.000014 & 0.000064 & 0.000001 & 0.000004 & 0.000001 & 0.000005 \\
\midrule
$m_3$ & \phantom{1}28 & 0.000215 & 0.000746 & 0.000215 & 0.000743 & 0.000008 & 0.000042 & 0.000011 & 0.000055 \\
$m_3$ & \phantom{1}55 & 0.000133 & 0.000523 & 0.000133 & 0.000519 & 0.000005 & 0.000028 & 0.000006 & 0.000035 \\
$m_3$ & 105 & 0.000090 & 0.000324 & 0.000090 & 0.000321 & 0.000003 & 0.000016 & 0.000004 & 0.000024 \\
$m_3$ & 210 & 0.000058 & 0.000191 & 0.000058 & 0.000192 & 0.000002 & 0.000012 & 0.000002 & 0.000016 \\
\midrule
$m_4$ & \phantom{1}28 & 0.000598 & 0.001452 & 0.000592 & 0.001461 & 0.000052 & 0.000122 & 0.000076 & 0.000162 \\
$m_4$ & \phantom{1}55 & 0.000351 & 0.000899 & 0.000349 & 0.000905 & 0.000030 & 0.000093 & 0.000041 & 0.000124 \\
$m_4$ & 105 & 0.000239 & 0.000520 & 0.000238 & 0.000516 & 0.000019 & 0.000052 & 0.000026 & 0.000064 \\
$m_4$ & 210 & 0.000149 & 0.000287 & 0.000147 & 0.000288 & 0.000014 & 0.000031 & 0.000020 & 0.000042 \\
\midrule
$m_5$ & \phantom{1}28 & 0.000862 & 0.003172 & 0.000861 & 0.003164 & 0.000040 & 0.000212 & 0.000050 & 0.000288 \\
$m_5$ & \phantom{1}55 & 0.000590 & 0.001977 & 0.000589 & 0.001974 & 0.000026 & 0.000116 & 0.000035 & 0.000132 \\
$m_5$ & 105 & 0.000393 & 0.001091 & 0.000396 & 0.001075 & 0.000017 & 0.000073 & 0.000022 & 0.000102 \\
$m_5$ & 210 & 0.000258 & 0.000621 & 0.000258 & 0.000618 & 0.000011 & 0.000039 & 0.000014 & 0.000048 \\
\midrule
$m_6$ & \phantom{1}28 & 0.000391 & 0.001927 & 0.000388 & 0.001927 & 0.000017 & 0.000136 & 0.000022 & 0.000166 \\
$m_6$ & \phantom{1}55 & 0.000237 & 0.001203 & 0.000238 & 0.001199 & 0.000011 & 0.000080 & 0.000012 & 0.000113 \\
$m_6$ & 105 & 0.000146 & 0.000666 & 0.000146 & 0.000668 & 0.000007 & 0.000042 & 0.000008 & 0.000054 \\
$m_6$ & 210 & 0.000093 & 0.000387 & 0.000093 & 0.000384 & 0.000004 & 0.000023 & 0.000006 & 0.000035 \\
\bottomrule
\end{tabular}
\end{table}

To assess the performance of both estimators near the boundary, a similar simulation was conducted, restricting the squared errors in the $\smash{\widetilde{\mathrm{ISE}}}$ values to points close to the boundary. Specifically, consider the simplex with a buffer $\delta \in (0, \infty)$, defined as
\[
\mathcal{S}_d(\delta) = \big\{\bb{s}\in \mathcal{S}_d : s_1 \geq \delta, \ldots, s_d \geq \delta, 1 - \|\bb{s}\|_1 \geq \delta\big\},
\]
and the corresponding boundary $\smash{\widetilde{\mathrm{ISE}}}$ values, viz.
\[
\text{boundary } \widetilde{\mathrm{ISE}}_{\mathrm{method}}^{j,r} = \frac{1}{1000} \sum_{i=1}^{1000} |\hat{m}_{n,\hat{b}_{n,r}}^{(\mathrm{method})}(\bb{V}_i) - m_j (\bb{V}_i)|^2 / 2, \vspace{-1mm}
\]
where $\mathcal{S}_d\backslash \mathcal{S}_d(\delta)$ denotes the boundary region, $\bb{V}_1,\ldots,\bb{V}_n\stackrel{\mathrm{iid}}{\sim} \mathcal{U}\{\mathcal{S}_d\backslash \mathcal{S}_d (n^{-1/3}/5)\}$, and $n^{-1/3}$ is the asymptotically optimal bandwidth order according to Proposition~\ref{prop:2}. In each case $(m_j, n$), Table~\ref{tab:2} again shows that the new LL estimator with Dirichlet kernel has the smallest mean, median, standard deviation, and interquartile range. \vspace{-1.5mm}

\begin{table}[ht]
\caption{Comparison of the LL and NW methods near the boundary, using the fixed design points $\bb{x}_1,\ldots,\bb{x}_n$, in terms of the mean, median, standard deviation (SD), and interquartile range (IQR), of $100$ boundary $\smash{\widetilde{\mathrm{ISE}}}$ values for regression functions $m_1, \ldots, m_6$ and sample size $n \in \{28, 55, 105, 210\}$.}
\label{tab:2}

\bigskip
\renewcommand{\arraystretch}{1.05}
\setlength{\tabcolsep}{3.3pt}
\centering
\begin{tabular}{cc|cccccccc}
\toprule
\multirow{2}{*}{Function} & \multirow{2}{*}{$n$} & \multicolumn{2}{c}{Mean} & \multicolumn{2}{c}{Median} & \multicolumn{2}{c}{SD} & \multicolumn{2}{c}{IQR} \\
\cmidrule(lr){3-4} \cmidrule(lr){5-6} \cmidrule(lr){7-8} \cmidrule(lr){9-10}
 & & LL & NW & LL & NW & LL & NW & LL & NW \\
\midrule
$m_1$ & \phantom{1}28 & 0.000055 & 0.000246 & 0.000055 & 0.000246 & 0.000002 & 0.000010 & 0.000003 & 0.000012 \\
$m_1$ & \phantom{1}55 & 0.000039 & 0.000133 & 0.000039 & 0.000133 & 0.000002 & 0.000006 & 0.000003 & 0.000008 \\
$m_1$ & 105 & 0.000029 & 0.000076 & 0.000028 & 0.000076 & 0.000001 & 0.000004 & 0.000002 & 0.000005 \\
$m_1$ & 210 & 0.000021 & 0.000046 & 0.000021 & 0.000046 & 0.000001 & 0.000002 & 0.000001 & 0.000003 \\
\midrule
$m_2$ & \phantom{1}28 & 0.000065 & 0.000472 & 0.000065 & 0.000472 & 0.000003 & 0.000020 & 0.000003 & 0.000027 \\
$m_2$ & \phantom{1}55 & 0.000058 & 0.000326 & 0.000059 & 0.000325 & 0.000003 & 0.000015 & 0.000004 & 0.000019 \\
$m_2$ & 105 & 0.000043 & 0.000226 & 0.000043 & 0.000226 & 0.000002 & 0.000011 & 0.000002 & 0.000013 \\
$m_2$ & 210 & 0.000025 & 0.000133 & 0.000026 & 0.000133 & 0.000001 & 0.000007 & 0.000002 & 0.000010 \\
\midrule
$m_3$ & \phantom{1}28 & 0.000315 & 0.001395 & 0.000315 & 0.001398 & 0.000012 & 0.000058 & 0.000017 & 0.000077 \\
$m_3$ & \phantom{1}55 & 0.000198 & 0.000851 & 0.000198 & 0.000845 & 0.000008 & 0.000050 & 0.000011 & 0.000078 \\
$m_3$ & 105 & 0.000135 & 0.000572 & 0.000135 & 0.000570 & 0.000005 & 0.000031 & 0.000006 & 0.000045 \\
$m_3$ & 210 & 0.000090 & 0.000381 & 0.000090 & 0.000378 & 0.000004 & 0.000018 & 0.000006 & 0.000022 \\
\midrule
$m_4$ & \phantom{1}28 & 0.001259 & 0.003628 & 0.001259 & 0.003599 & 0.000069 & 0.000182 & 0.000095 & 0.000218 \\
$m_4$ & \phantom{1}55 & 0.000855 & 0.002365 & 0.000855 & 0.002364 & 0.000051 & 0.000129 & 0.000067 & 0.000153 \\
$m_4$ & 105 & 0.000639 & 0.001675 & 0.000634 & 0.001679 & 0.000041 & 0.000100 & 0.000052 & 0.000111 \\
$m_4$ & 210 & 0.000447 & 0.001121 & 0.000446 & 0.001121 & 0.000026 & 0.000068 & 0.000039 & 0.000090 \\
\midrule
$m_5$ & \phantom{1}28 & 0.001140 & 0.006616 & 0.001132 & 0.006599 & 0.000059 & 0.000370 & 0.000081 & 0.000472 \\
$m_5$ & \phantom{1}55 & 0.000790 & 0.003834 & 0.000792 & 0.003835 & 0.000038 & 0.000209 & 0.000047 & 0.000258 \\
$m_5$ & 105 & 0.000531 & 0.002315 & 0.000531 & 0.002304 & 0.000026 & 0.000122 & 0.000035 & 0.000164 \\
$m_5$ & 210 & 0.000349 & 0.001439 & 0.000348 & 0.001438 & 0.000015 & 0.000076 & 0.000023 & 0.000091 \\
\midrule
$m_6$ & \phantom{1}28 & 0.000731 & 0.004082 & 0.000730 & 0.004055 & 0.000025 & 0.000179 & 0.000035 & 0.000255 \\
$m_6$ & \phantom{1}55 & 0.000483 & 0.002377 & 0.000484 & 0.002348 & 0.000017 & 0.000123 & 0.000020 & 0.000159 \\
$m_6$ & 105 & 0.000318 & 0.001410 & 0.000319 & 0.001414 & 0.000012 & 0.000079 & 0.000018 & 0.000108 \\
$m_6$ & 210 & 0.000218 & 0.000887 & 0.000219 & 0.000883 & 0.000009 & 0.000050 & 0.000014 & 0.000057 \\
\bottomrule
\end{tabular}
\end{table}

To further illustrate the superior performance of the LL smoother, an additional simulation was performed using a random sample of design points $\bb{X}_1, \ldots, \bb{X}_n$ drawn from a $\mathrm{Dirichlet}(2 \times \bb{1}, 2)$ distribution, instead of the fixed design points $\bb{x}_1,\ldots,\bb{x}_n$ from the mesh $\mathcal{M}_k$. The $\smash{\widetilde{\mathrm{ISE}}}$ values were adjusted to account for the random sampling by incorporating density weights, viz.
\[
\text{weighted } \widetilde{\mathrm{ISE}}_{\mathrm{method}}^{j,r} = \frac{1}{1000} \sum_{i=1}^{1000} |\hat{m}_{n,\hat{b}_{n,r}}^{(\mathrm{method})}(\bb{U}_i) - m_j(\bb{U}_i)|^2 f(\bb{U}_i) / 2.
\]
Table~\ref{tab:3} displays the results in each case $(m_j,n$), showing once more that the new LL estimator with Dirichlet kernel has the smallest mean, median, standard deviation, and interquartile range.

\begin{table}[ht]
\caption{Comparison of the LL and NW methods, using a random sample of design points $\bb{X}_1,\ldots,\bb{X}_n$ from a $\mathrm{Dirichlet}(2 \times \bb{1}, 2)$ distribution, in terms of the mean, median, standard deviation (SD), and interquartile range (IQR), of $100$ weighted $\smash{\widetilde{\mathrm{ISE}}}$ values for regression functions $m_1, \ldots, m_6$ and sample size $n \in \{28, 55, 105, 210\}$.}
\label{tab:3}

\bigskip
\renewcommand{\arraystretch}{1.05} 
\setlength{\tabcolsep}{3.3pt} 
\centering
\begin{tabular}{cc|cccccccc}
\toprule
\multirow{2}{*}{Function} & \multirow{2}{*}{$n$} & \multicolumn{2}{c}{Mean} & \multicolumn{2}{c}{Median} & \multicolumn{2}{c}{SD} & \multicolumn{2}{c}{IQR} \\
\cmidrule(lr){3-4} \cmidrule(lr){5-6} \cmidrule(lr){7-8} \cmidrule(lr){9-10}
 & & LL & NW & LL & NW & LL & NW & LL & NW \\
\midrule
$m_1$ & \phantom{1}28 & 0.000124 & 0.000473 & 0.000111 & 0.000440 & 0.000052 & 0.000160 & 0.000073 & 0.000200 \\
$m_1$ & \phantom{1}55 & 0.000059 & 0.000341 & 0.000054 & 0.000335 & 0.000023 & 0.000091 & 0.000018 & 0.000113 \\
$m_1$ & 105 & 0.000035 & 0.000293 & 0.000033 & 0.000286 & 0.000008 & 0.000061 & 0.000009 & 0.000071 \\
$m_1$ & 210 & 0.000024 & 0.000263 & 0.000024 & 0.000259 & 0.000004 & 0.000041 & 0.000005 & 0.000052 \\
\midrule
$m_2$ & \phantom{1}28 & 0.000058 & 0.001373 & 0.000056 & 0.001289 & 0.000017 & 0.000550 & 0.000029 & 0.000593 \\
$m_2$ & \phantom{1}55 & 0.000033 & 0.001079 & 0.000032 & 0.001055 & 0.000007 & 0.000300 & 0.000010 & 0.000339 \\
$m_2$ & 105 & 0.000021 & 0.000838 & 0.000021 & 0.000803 & 0.000004 & 0.000222 & 0.000005 & 0.000262 \\
$m_2$ & 210 & 0.000013 & 0.000780 & 0.000013 & 0.000758 & 0.000002 & 0.000139 & 0.000003 & 0.000154 \\
\midrule
$m_3$ & \phantom{1}28 & 0.000279 & 0.004446 & 0.000266 & 0.004051 & 0.000076 & 0.001643 & 0.000108 & 0.001974 \\
$m_3$ & \phantom{1}55 & 0.000168 & 0.003263 & 0.000161 & 0.003186 & 0.000039 & 0.000832 & 0.000045 & 0.001039 \\
$m_3$ & 105 & 0.000102 & 0.002869 & 0.000099 & 0.002831 & 0.000017 & 0.000680 & 0.000022 & 0.000957 \\
$m_3$ & 210 & 0.000066 & 0.002512 & 0.000065 & 0.002478 & 0.000008 & 0.000399 & 0.000012 & 0.000589 \\
\midrule
$m_4$ & \phantom{1}28 & 0.000616 & 0.003734 & 0.000582 & 0.003471 & 0.000163 & 0.001358 & 0.000207 & 0.001407 \\
$m_4$ & \phantom{1}55 & 0.000362 & 0.002856 & 0.000350 & 0.002737 & 0.000071 & 0.000918 & 0.000090 & 0.001120 \\
$m_4$ & 105 & 0.000227 & 0.002505 & 0.000222 & 0.002484 & 0.000039 & 0.000663 & 0.000046 & 0.000679 \\
$m_4$ & 210 & 0.000165 & 0.002242 & 0.000162 & 0.002210 & 0.000020 & 0.000351 & 0.000026 & 0.000457 \\
\midrule
$m_5$ & \phantom{1}28 & 0.001232 & 0.014832 & 0.001156 & 0.013689 & 0.000360 & 0.005644 & 0.000507 & 0.006522 \\
$m_5$ & \phantom{1}55 & 0.000748 & 0.010920 & 0.000722 & 0.010300 & 0.000187 & 0.003242 & 0.000244 & 0.004549 \\
$m_5$ & 105 & 0.000481 & 0.009099 & 0.000474 & 0.008865 & 0.000076 & 0.002257 & 0.000100 & 0.002248 \\
$m_5$ & 210 & 0.000298 & 0.008303 & 0.000294 & 0.008358 & 0.000032 & 0.001360 & 0.000048 & 0.001716 \\
\midrule
$m_6$ & \phantom{1}28 & 0.000596 & 0.008509 & 0.000575 & 0.007737 & 0.000218 & 0.003111 & 0.000283 & 0.003400 \\
$m_6$ & \phantom{1}55 & 0.000327 & 0.006492 & 0.000310 & 0.006237 & 0.000085 & 0.001949 & 0.000094 & 0.002822 \\
$m_6$ & 105 & 0.000189 & 0.005341 & 0.000187 & 0.005263 & 0.000035 & 0.001279 & 0.000049 & 0.001731 \\
$m_6$ & 210 & 0.000112 & 0.005060 & 0.000108 & 0.005025 & 0.000016 & 0.000829 & 0.000022 & 0.001155 \\
\bottomrule
\end{tabular}
\end{table}

In their paper, \citet{MR4796622} also considered a NW regression estimator where the Dirichlet kernel is replaced by a discrete multinomial kernel, offering a Bernstein smoothing alternative. No advantages were found over the NW estimator with Dirichlet kernel \eqref{eq:4}, so it was not included in the present study for the sake of concision. The LL method with Dirichlet kernel outperforms both NW estimators by a significant margin. This is unsurprising, however, given that NW estimators can be viewed as local {\it constant} estimators, where the slope of the corresponding locally weighted least-squares problem is set to zero.

\section{Data example}
\label{sec:6}

As a concrete illustration, the proposed local linear (LL) smoother  is used to take yet another look at a classical dataset reported by \citet[p.~359]{MR865647} and displayed in Table~\ref{tab:Aitchison.data.5} for convenience.

\begin{table}[b!]
\centering
\caption{Sand, silt, clay compositions of 39 sediment samples at different water depths in Stanwell-Fletcher Lake, Somerset Island (Inuktitut Kuuganajuk, Nunavut, Canada).}
\label{tab:Aitchison.data.5} 

\bigskip
\renewcommand{\arraystretch}{1.05} 
\setlength{\tabcolsep}{3.3pt} 
\centering
\begin{tabular}{ccccc|ccccc}
\toprule
Sediment & \multicolumn{3}{c}{Percentages} & Water depth & Sediment & \multicolumn{3}{c}{Percentages} & Water depth \\
no. & Sand & Silt & Clay & (meters) & no. & Sand & Silt & Clay & (meters) \\
\midrule
\phantom{0}1 & 77.5 & 19.5 & \phantom{0}3.0 & 10.4 & 21 & \phantom{0}9.5 & 53.5 & 37.0 & \phantom{0}47.1 \\
\phantom{0}2 & 71.9 & 24.9 & \phantom{0}3.2 & 11.7 & 22 & 17.1 & 48.0 & 34.9 & \phantom{0}48.4 \\
\phantom{0}3 & 50.7 & 36.1 & 13.2 & 12.8 & 23 & 10.5 & 55.4 & 34.1 & \phantom{0}49.4 \\
\phantom{0}4 & 52.2 & 40.9 & \phantom{0}6.6 & 13.0 & 24 & \phantom{0}4.8 & 54.7 & 41.0 & \phantom{0}49.5 \\
\phantom{0}5 & 70.0 & 26.5 & \phantom{0}3.5 & 15.7 & 25 & \phantom{0}2.6 & 45.2 & 52.2 & \phantom{0}59.2 \\
\phantom{0}6 & 66.5 & 32.2 & \phantom{0}1.3 & 16.3 & 26 & 11.4 & 52.7 & 35.9 & \phantom{0}60.1 \\
\phantom{0}7 & 43.1 & 55.3 & \phantom{0}1.6 & 18.0 & 27 & \phantom{0}6.7 & 46.9 & 46.4 & \phantom{0}61.7 \\
\phantom{0}8 & 53.4 & 36.8 & \phantom{0}9.8 & 18.7 & 28 & \phantom{0}6.9 & 49.7 & 43.4 & \phantom{0}62.4 \\
\phantom{0}9 & 15.5 & 54.4 & 30.1 & 20.7 & 29 & \phantom{0}4.0 & 44.9 & 51.1 & \phantom{0}69.3 \\
10 & 31.7 & 41.5 & 26.8 & 22.1 & 30 & \phantom{0}7.4 & 51.6 & 40.9 & \phantom{0}73.6 \\
11 & 65.7 & 27.8 & \phantom{0}6.5 & 22.4 & 31 & \phantom{0}4.8 & 49.5 & 45.7 & \phantom{0}74.4 \\
12 & 70.4 & 29.0 & \phantom{0}0.6 & 24.4 & 32 & \phantom{0}4.5 & 48.5 & 47.0 & \phantom{0}78.5 \\
13 & 17.4 & 53.6 & 29.0 & 25.8 & 33 & \phantom{0}6.6 & 52.1 & 41.3 & \phantom{0}82.9 \\
14 & 10.6 & 69.8 & 19.6 & 32.5 & 34 & \phantom{0}6.7 & 47.3 & 45.9 & \phantom{0}87.7 \\
15 & 38.2 & 43.1 & 18.7 & 33.6 & 35 & \phantom{0}7.4 & 45.6 & 46.9 & \phantom{0}88.1 \\
16 & 10.8 & 52.7 & 36.5 & 36.8 & 36 & \phantom{0}6.0 & 48.9 & 45.1 & \phantom{0}90.4 \\
17 & 18.4 & 50.7 & 30.9 & 37.8 & 37 & \phantom{0}6.3 & 53.8 & 39.9 & \phantom{0}90.6 \\
18 & \phantom{0}4.6 & 47.4 & 48.0 & 36.9 & 38 & \phantom{0}2.5 & 48.0 & 49.5 & \phantom{0}97.7 \\
19 & 15.6 & 50.4 & 34.0 & 42.2 & 39 & \phantom{0}2.0 & 47.8 & 50.2 & 103.7 \\
20 & 31.9 & 45.1 & 23.0 & 47.0 & & & & & \\
\bottomrule
\end{tabular}
\end{table}

These data, which are also available in the \textsf{R} package \texttt{DirichletReg}, were collected by \cite{Coakler/Rust:1968}, who provide the composition in terms of sand, silt, and clay of 39 samples of sediment in Stanwell-Fletcher Lake, Somerset Island (Inuktitut Kuuganajuk, Nunavut, Canada) as a function of depth (in meters).

One way to assess the relation between water depth and compositional pattern is to  determine the extent to which the former can be predicted by the latter, which calls for a model of the form~\eqref{eq:1}.

The dataset comprises 39 design points $\bb{x}_i = (x_{i,1}, x_{i,2})$, which represent the proportions of sand and silt in each sediment sample, respectively. The proportion of clay is determined by the complement, $1 - x_{i,1} - x_{i,2}$. Water depth (in meters) of each sediment sample, denoted $y_i$, is treated as the dependent variable.

\begin{figure}[!t]
\centering
\includegraphics[width=0.90\textwidth]{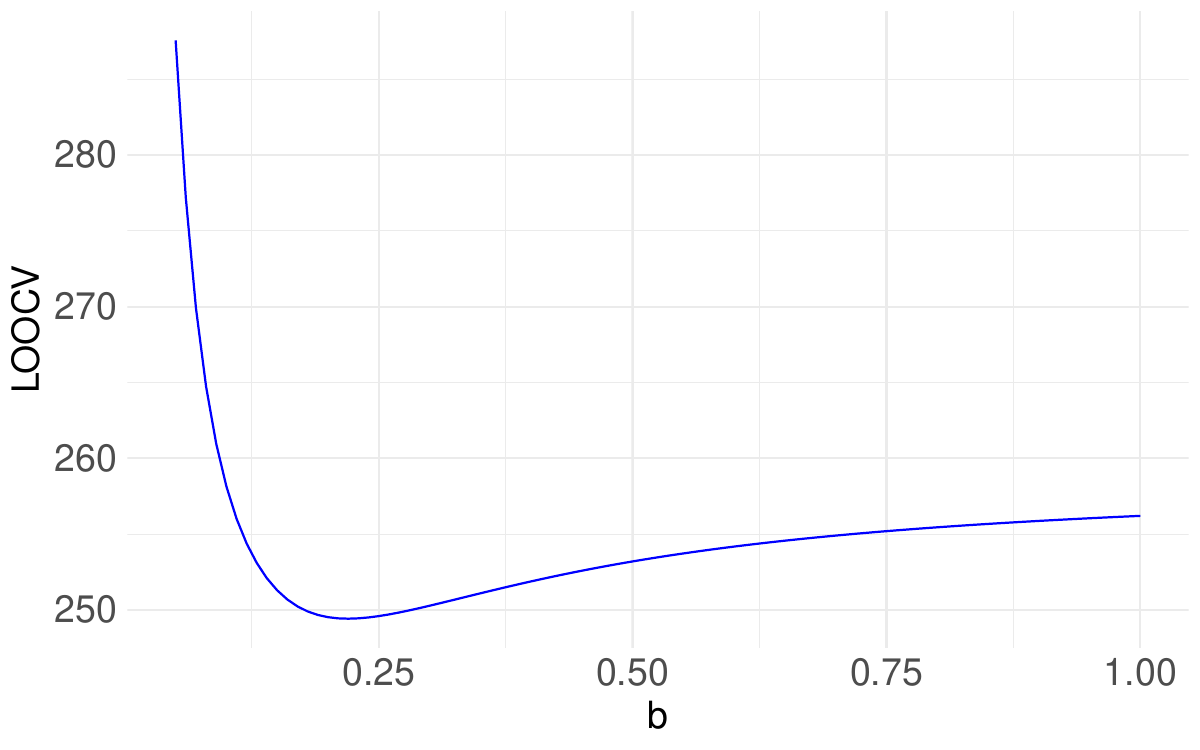}
\caption{Plot of the leave-one-out cross-validation as a function of the bandwidth $b$.}
\label{fig:LOOCV.b}
\end{figure}

The local linear smoother $\smash{\hat{m}_{n,b}^{(\mathrm{LL})}}$, as defined in Section~\ref{sec:2}, is employed to estimate the water depth based of sediment composition. The bandwidth $b$ is selected using the leave-one-out cross-validation method which minimizes
\[
\mathrm{LOOCV}(b) = \frac{1}{n} \sum_{i=1}^n \{y_i - \hat{m}_{n,b,(-i)}^{(\mathrm{LL})}(\bb{x}_i)\}^2,
\]
where $\smash{\hat{m}_{n,b,(-i)}^{(\mathrm{LL})}}$ represents the leave-one-out estimator for each $i\in \{1,\ldots,n\}$, i.e., the local linear smoother defined without the $i$th pair $(\bb{x}_i, y_i)$. Fig.~\ref{fig:LOOCV.b} shows the graph of $\mathrm{LOOCV}$ as a function of the bandwidth $b$.

\begin{figure}[!t]
\centering
\includegraphics[width=0.85\textwidth]{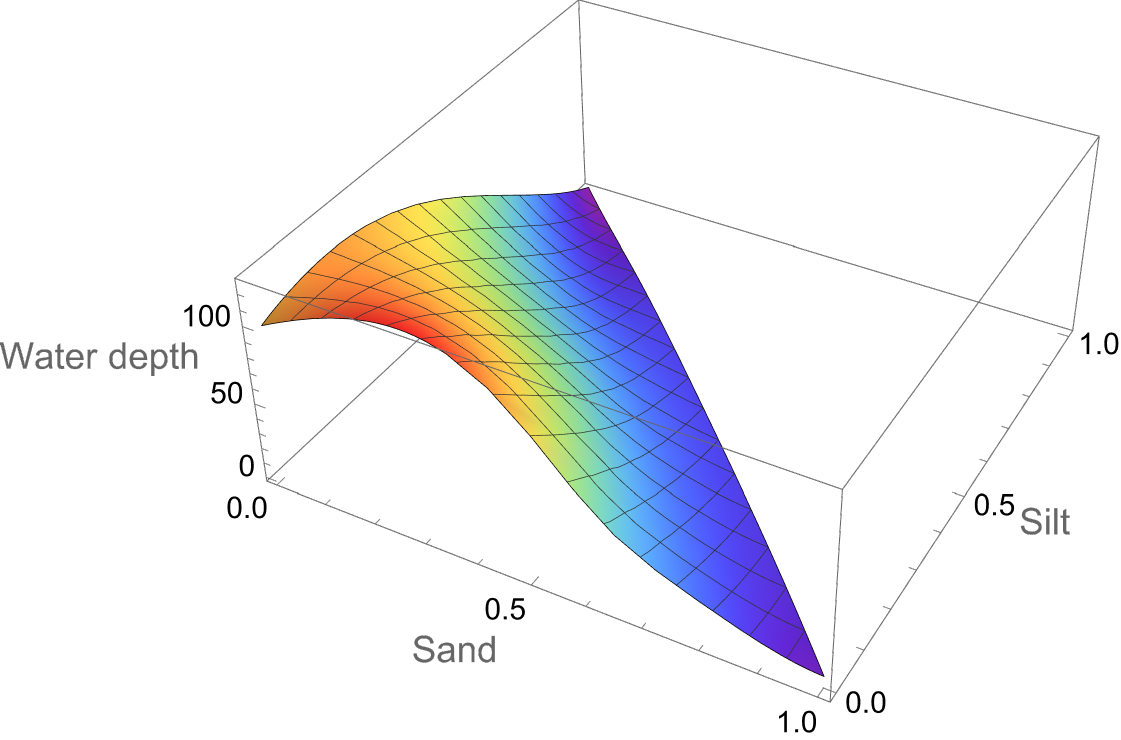}
\caption{Three-dimensional plot of estimated water depth as a function of sand and silt proportions.}
\label{fig:3d.plot}
\end{figure}

\begin{figure}[!b]
\centering
\includegraphics[width=0.85\textwidth]{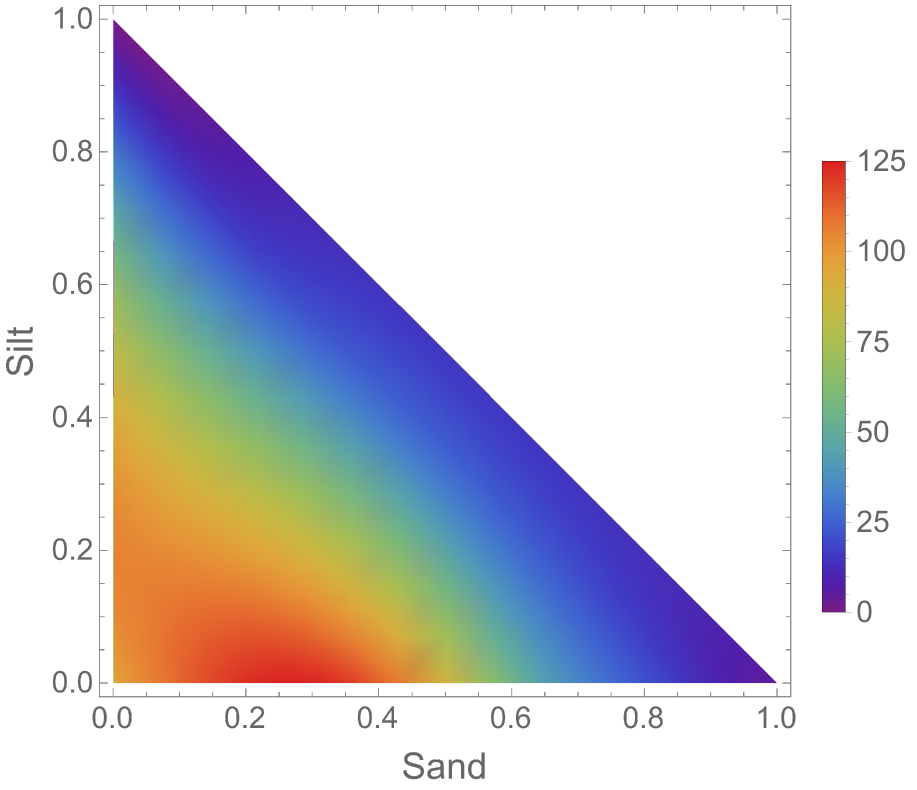}
\caption{Density plot of the estimated water depth as a function of the proportion of sand and silt.}
\label{fig:density.plot}
\end{figure}

Numerical computations reveal that the optimal bandwidth under this criterion is
\[
\hat{b} = \operatorname{argmin}_{b\in (0,\infty)} \mathrm{LOOCV}(b) \approx 0.2195. 
\]
The 3-dimensional plot and density plot of $\smash{\hat{m}_{39,\hat{b}}^{(\mathrm{LL})}}$ are illustrated in Fig.~\ref{fig:3d.plot}~and~Fig.~\ref{fig:density.plot}, respectively. These graphs make it clear that sand and silt are present in smaller proportions than clay in sediment composition in deeper areas of the lakebed. In the deepest region (in red), sediments are predominantly clay, with the remaining composition skewing heavily toward sand.

\section{Conclusion and outlook}
\label{sec:7}

This paper introduced and studied a local linear smoother for nonparametric regression on the simplex using Dirichlet kernels. The proposed method leverages asymmetric kernels and local linear weighting to mitigate boundary bias, providing a natural multivariate extension of Chen's local linear smoother with beta kernel \citep{MR1910175}. The estimator is versatile and can be applied in both fixed and random design settings, making it suitable for a wide range of practical scenarios. The theoretical findings include asymptotic expressions for the bias, variance, and mean squared error of the estimator, derived both at the center of the simplex and near its boundary, as well as an expression for the mean integrated squared error, which generalizes Chen's work to regression surfaces on the simplex. The simulations in Section~\ref{sec:5} show convincing improved performance compared to the existing Nadaraya--Watson estimator with Dirichlet kernels proposed by \citet{MR4796622}.

Looking ahead, there are several interesting directions for future research. First, the local linear approach could be extended to local polynomial estimation of higher degree, potentially achieving greater bias reduction and more accurate local approximations; refer to Section~4 of \citet{MR1311979} for seminal results in this direction in the classical local polynomial setting.

Second, establishing the uniform strong consistency of $\smash{\hat{m}^{(\mathrm{LL})}_{n,b}}$ over a growing sequence of compact sets inside the simplex would be valuable for both theoretical innovation and practical guidance. In particular, such rates would clarify how the estimator behaves as the dimension $d$ grows, and could provide insights for designing adaptive bandwidth selection procedures in a minimax sense \citep{MR1705298,MR2724359}. For reference, in the context of the asymmetric kernel density estimator with Dirichlet kernel, the uniform strong consistency on growing compacts was shown in Theorem~4 of \citet{MR4319409} using a chaining argument, but the rate of convergence obtained there is not optimal.

Third, deriving the asymptotic distribution of the local linear smoother would enable the construction of confidence intervals for $m(\bb{s})$ and facilitate hypothesis testing and inference.

Fourth, extending both the theoretical results established in the present paper and the open problems mentioned above --- including uniform strong consistency, asymptotic normality, and other related properties --- to the case of dependent observations would allow for a more accurate analysis of data with spatial dependence, such as neighboring soil samples, or time dependence, such as time series tracking the evolution of a compositional process and the corresponding predictions.

All four directions mentioned above remain open problems and provide fertile ground for future research. From a computational perspective, ensuring the scalability and efficiency of the proposed method is critical, especially for large datasets and higher-dimensional simplices. For example, the bandwidth selection procedure can be made faster in higher dimensions using parallel computing, dimension reduction techniques, using a pilot estimate for $h_{\emptyset}$ to compute a plug-in version of the globally optimal bandwidth $b_{\mathrm{opt}}$ in Proposition~\ref{prop:2}, and so forth.

By pursuing these extensions, the aim is to enrich the theoretical foundations of Dirichlet kernel-based local linear smoothing and to broaden its range of applications.

\section*{Appendix: Technical results}

The first lemma recalls a well-known formula from linear algebra for the inverse of a $2\times 2$ block matrix, as detailed, for example, in Section~0.7 of \citet{MR2978290}.

\bigskip
\begin{lemma}\label{lem:inverse.block.matrix}
Let $M$ be a square matrix of width at least $2$ which is partitioned into a $2\times 2$ block matrix. If the upper left block $M_{11}$ is invertible and the Schur complement $M_{\star} = M_{22} - M_{21} M_{11}^{-1} M_{12}$ is invertible, then
\[
M =
\begin{bmatrix}
M_{11} & ~ M_{12} \\[1mm]
M_{21} & ~ M_{22}
\end{bmatrix}
~ \Rightarrow ~~
M^{-1} =
\begin{bmatrix}
M_{11}^{-1} + M_{11}^{-1} M_{12} M_{\star}^{-1} M_{21} M_{11}^{-1} & ~ - M_{11}^{-1} M_{12} M_{\star}^{-1} \\[1mm]
- M_{\star}^{-1} M_{21} M_{11}^{-1} & ~ M_{\star}^{-1}
\end{bmatrix}\!.
\]
\end{lemma}

The second lemma, which corresponds to Lemma~1 of \citet{MR4319409}, provides a uniform bound on $A_b(\bb{s}) = \EE\{\kappa_{\bb{s},b}(\bb{X})^2\}$ and describes its asymptotic behavior as $b\to 0$, whether the sequence of points $\bb{s} = \bb{s}(b)$ is well inside or near the boundary of the simplex $\mathcal{S}_d$.

\bigskip
\begin{lemma}\label{lem:A.b.asymptotics}
Uniformly for $\bb{s}\in \mathcal{S}_d$, one has, as $b\to 0$,
\[
0 < A_b(\bb{s}) \leq \frac{b^{(d + 1) / 2} \, (1 / b + d)^{d + 1/2}}{(4\pi)^{d/2} \sqrt{(1 - \|\bb{s}\|_1) \prod_{i\in [d]} s_i}} \, \{1 + \OO(b)\}.
\]
Also, for any vector $\bb{\lambda} \in (0, \infty)^d$ and any subset $\mathcal{J}\subseteq [d]$ of indices and corresponding map $\psi_{\mathcal{J}}$ defined in \eqref{eq:3}, one has, as $b\to 0$,
\[
A_b(\bb{s}) = b^{-d/2} \, \psi_{\emptyset}(\bb{s}) \times \{1 + \OO_{\bb{s}}(b)\},
\]
if $s_i$ is fixed for all $i \in [d]$ and $(1 - \|\bb{s}\|_1)$ is fixed, while
\[
A_b(\bb{s}) = b^{-(d + |\mathcal{J}|)/2} \psi_{\mathcal{J}}(\bb{s}) \prod_{i\in \mathcal{J}} \frac{\Gamma(2\lambda_i + 1)}{2^{2\lambda_i + 1} \Gamma^2(\lambda_i + 1)} \times \{1 + \OO_{\bb{\lambda},\bb{s}}(b)\}
\]
if $s_i / b \to \lambda_i \text{ for all } i \in \mathcal{J}$, $s_i$ is fixed for all $i \in [d]\backslash \mathcal{J}$ and $(1 - \|\bb{s}\|_1)$ is fixed.
\end{lemma}

\section*{Reproducibility}\label{sec:reproducibility}
\addcontentsline{toc}{section}{Reproducibility}

The \textsf{R} codes used to produce the figures, conduct the simulations, and analyze the data example can be accessed in the GitHub repository of \citet{GenestOuimet2024github}.

\section*{Acknowledgments}
\addcontentsline{toc}{section}{Acknowledgments}

The simulations for this study were conducted using the computational resources of Calcul Québec (\href{https://www.calculquebec.ca}{www.calculquebec.ca}) and the Digital Research Alliance of Canada (\href{https://www.alliancecan.ca}{www.alliancecan.ca}). The authors are grateful to the four reviewers for their comprehensive and insightful comments, which enhanced the quality of this manuscript.

\section*{Funding}
\addcontentsline{toc}{section}{Funding}

Genest's research was funded in part by the Canada Research Chairs Program (Grant no.~950--231937) and the Natural Sciences and Engineering Research Council of Canada (RGPIN--2024--04088). Ouimet's funding was made possible through a contribution to Genest's work from the Trottier Institute for Science and Public Policy.

\phantomsection
\addcontentsline{toc}{chapter}{References}

{\small
\bibliographystyle{authordate1}
\bibliography{bib}
}

\end{document}